\documentclass[12pt]{amsart}
\usepackage{graphicx}
\usepackage[all]{xy}
\usepackage{epstopdf}
\newtheorem{theorem}{Theorem}[section]
\newtheorem{corollary}[theorem]{Corollary}

\newtheorem{claim}[theorem]{Claim}

\newtheorem{proposition}[theorem]{Proposition}
\theoremstyle{definition}
\newtheorem{definition}[theorem]{Definition}

\theoremstyle{remark}

\newtheorem{example}[theorem]{Example}
\numberwithin{equation}{section}
\makeatletter \def\qtri{\mathinner{{\bf Q}\mkern-12mu\raise2.5\p@ \vbox{\kern1\p@\hbox{{\tiny{$\rhd$}} }}}}\makeatother
\begin{document}
\bibliographystyle{plain}

\title{Toward a Classification of Finite Quandles}
\author{Geoff Ehrman, Ata Gurpinar\\Matthew Thibault, David Yetter}
\thanks{The  authors were partially supported by the Kansas State
University REU and NSF grant GOMT530725.\\
\\Department of Mathematics, Kansas State University \\
Manhattan, Kansas 66506, USA}
\maketitle

\noindent{\small {\bf Abstract:} {\em This paper summarizes substantive new results derived by a student team (the first three authors) under the direction of the fourth author at the 2005 session of the KSU REU ``Brainstorming and Barnstorming''.  The main results are a decomposition theorem for quandles in terms of
an operation of `semidisjoint union' showing that all finite quandles canonically decompose via iterated semidisjoint unions into connected subquandles, and a structure theorem for finite connected quandles with prescribe inner automorphism group.
The latter theorem suggests a new approach to the classification of finite connected quandles.}}

\section{Introduction} 

Quandles were introduced by Joyce \cite{Joyce1} \cite{Joyce2} as an
algebraic invariant of classical knots and links.

\begin{definition} \cite{Joyce1} \cite{Joyce2}
A {\em quandle} $(Q,\rhd)$ is a set $Q$ equipped with an binary operation $\rhd$ such that the following conditions hold:
\begin{itemize}
\item every element in the set is {\em idempotent} with respect to $\rhd$: $\forall x \in Q, x \rhd x = x$
\item $\rhd$ is {\em invertible} as a right-acting operator: $\forall a, b \in Q, \exists ! y \in Q : y \rhd a = b$
\item $\rhd$ is {\em right-distributive} over itself: $\forall a,b,c \in Q, (a \rhd b) \rhd c = (a \rhd c) \rhd (b \rhd c)$
\end{itemize}
\end{definition}

When the operation is clear from context, we will denote the quandle by its underlying set.

The second axiom is equivalent to the existence of a second operation $\rhd^{-1}$ for which $\rhd^{-1} y$
is the inverse to $\rhd y$ for all $y \in Q$.

Quandles may be regarded as an abstraction from groups in as much as many important examples arise as conjugation invariant subsets of groups with the operation $\rhd y$ given by right conjugation by $y$.  (Indeed it is a theorem of Joyce \cite{Joyce1}
that free quandles are isomorphic to disjoint unions of conjugacy classes in groups, and thus that the equational theory of quandles is `the equational theory of groups under conjugation'.) 

The knot quandle \cite{Joyce1} \cite{Joyce2}, though it admits a homotopy theoretic definition, can be described most simply by generators
and relations:  modify the Wirtinger presentation of the
knot group by replacing right conjugation with $\rhd$ and
left conjugation with $\rhd^{-1}$.

Finite quandles, in particular, are of some interest, since
they give rise both to `counting invariants', of classical knots
and link, which generalize
Tait's notion of three-coloring a knot, to `'counting invariants' of monodromy situations (cf. \cite{Yetter.monodromy}), and form the basis for more refined
topological invariants derived from quandle cohomology
(cf. \cite{CJKLS, CJKS}).

The present work is intended as a contribution to the problem
of classifying finite quandles.

\section{Notation, Examples, and Basic Concepts}

\begin{definition}
A {\em quandle homomorphism}, given $(Q, \rhd)$ and $(Q', \rhd')$, is a mapping $\rho: Q \rightarrow Q'$ such that $\forall x, y \in Q, (x \rhd y)\rho =(x)\rho \rhd' (y)\rho$.

Note, for the sake of agreement with the action of elements of a quandle
on the quandle itself, which is
written as a right action by quandle homomorphisms, we write quandle 
homomorphisms to the right of their arguments.
  
\end{definition}
 
\begin{example}
${\bf n}$ is the trivial quandle of order n. Given $|Q|=n, {\bf n} :=(Q, \rhd)$ where $\forall x, y \in Q, x \rhd y = x$.
\end{example}
\begin{example}
The Tait quandle $({\bf T}_3, \rhd)$ is the quandle with underlying
set $\{a, b, c\}$ and operation

\smallskip
\begin{center}
\begin{tabular}{c|ccc}
$\rhd$ & $a$ & $b$ & $c$ \\ \hline
$a$ & $a$ & $c$ & $b$ \\
$b$ & $c$ & $b$ & $a$ \\
$c$ & $b$ & $a$ & $c$ 
\end{tabular}
\end{center}
\smallskip

This quandle is so named because Tait's notion of three-coloring a knot is 
equivalent to the existence of a non-trivial quandle homomorphism from the
knot quandle (cf. Joyce \cite{Joyce1}) to ${\bf T}_3$.
\end{example}

\begin{definition}
Quandles $(Q, \rhd)$ and $(Q', \rhd')$ are {\em isomorphic} when there exists a bijective quandle homomorphism $\rho$ (the {\em isomorphism} from $Q$ to $Q'$. We denote the existence of such an isomorphism by $(Q, \rhd) \cong (Q', \rhd')$.
\end{definition}

\begin{definition}
The {\em automorphism group} of a quandle $(Q, \rhd)$, denoted {\em $Aut(Q)$}, is the group of all isomorphisms $\rho: Q \rightarrow Q$. The elements of $Aut(Q)$ act on those of $Q$ by right action.
\end{definition}

\begin{definition}
The {\em inner automorphism group} of a quandle $(Q, \rhd)$, denoted {\em $Inn(Q)$}, is the subgroup of $Aut(Q)$ generated by all $S_x$, where $\forall x,y \in Q, S_x(y):= y \rhd x$.
\end{definition}

\begin{definition}
The {\em orbit} of $s \in Q$ is the subset of elements $t \in Q$ such that there exists some $p \in Inn(Q)$ where $p$ maps $s$ to $t$.
\end{definition}

\begin{definition}
A quandle $(Q,\rhd)$ is {\em (algebraically) connected} when there exists exactly one orbit in $Q$--that is, $\forall x \in Q$, the orbit of $x$ {\em is} all of $Q$. 
\end{definition}

\begin{definition}
Given a set $Q$ and a group $G$ with a right action by quandle homomorphisms on $Q$, an {\em augmentation map} is a map $| \cdot |: Q \rightarrow G$ such that the following hold:
\begin{itemize}
\item $\forall q \in Q, q |q| = q$
\item $\forall q \in Q, g \in G, |qg| = g^{-1} |q| g$
\end{itemize}
\end{definition}

\begin{definition}
 The {\em universal augmentation group} of the quandle $Q$, denoted, 
$\Gamma_Q$, is the group freely generated by all formal augmentations $|x|$ 
of the elements $x \in Q$ 
modulo relations $\forall x, y \in Q, |x \rhd y| = |y|^{-1}|x||y|$.
\end{definition}

Joyce \cite{Joyce1} showed that the inclusion of generators into $\Gamma_Q$ is
the universal augmentation of $Q$ in the sense that the quandle operation
induces an action of
$\Gamma_Q$ on $Q$ by quandle homomorphisms such that the inclusion of
generators is an augmentation, and given any other
augmentation $\langle \cdot \rangle:Q\rightarrow G$, there is a unique group homomorphism
$c:\Gamma_Q\rightarrow G$ such that 
$\forall q \in Q, c(|q|) = \langle q \rangle$ and
$\forall q \in Q, \forall g \in \Gamma_Q, q c(g) = qg$.

\section{Semidisjoint Union of Quandles}
Let $(Q_1,\rhd_1)$, $(Q_2,\rhd_2), \ldots, (Q_n,\rhd_n)$
be quandles. For each quandle $(Q_i, \rhd_i)$, we have the universal 
augmentation map
$|\cdot|$ from it to $\Gamma_{Q_i}$.  In particular, for all $x, y \in
Q_i$, $x |y| := x \rhd_i y$. Note that this is an augmentation map
since $\forall x \in Q_i, x |x| = x \rhd_i x = x$ and, given some $g = |y_1|^{\pm 1} |y_2|^{\pm 1} \ldots |y_k|^{\pm 1} \in \Gamma_{Q_i}$,\\

\begin{align*}
|x g|= |x |y_1|^{\pm 1}|y_2|^{\pm 1} \ldots |y_k|^{\pm 1}| & =
|(\ldots(x \rhd_i^{\pm 1} y_1) \rhd_i^{\pm 1} y_2) \ldots
\rhd_i^{\pm 1} y_k)|
\\ & = |y_k|^{\mp 1} \ldots |y_1|^{\mp 1} |x| |y_1|^{\pm 1} \ldots
|y_k|^{\pm 1} 
\\ & = g^{-1} |x| g.
\end{align*}

Now, observe that the orbits of $Inn(Q)$ are subquandles of any quandle $Q$, and that if $Q$ has orbits $Q_1, \ldots ,Q_n$, 
the augmentation of $Q$ in $Inn(Q)$ induces an augmentation 
of each $Q_i$ in $Inn(Q)$, and thus a group homomorphism from $\Gamma_{Q_i}$ to
$Inn(Q)$.  These, in turn induce group homomorphisms from $\Gamma_{Q_i}$ to
$Aut(Q_j)$ for each $j$ (with the homomorphism from $\Gamma_{Q_i}$ to 
$Aut(Q_i)$ being that induced by the universal property of $\Gamma_{Q_i}$).
Notice, for $i \neq j$, these do not necessarily factor through the subgroup 
$Inn(Q_j)$.

This observation suggests the following construction:

For $1 \leq i, j \leq n$, let $g_{i,j}$ be a group homomorphism from $\Gamma_{Q_i}$ to $Aut(Q_j)$, such that $g_{i,i}$ is the canonical group 
homomorphism from
$\Gamma_{Q_i}$ to $Aut(Q_i)$.  Let $G$ be the $n \times n$ matrix of group
homomorphisms
with entries $g_{i,j}$. Define the operation $$\#(Q_1, Q_2, \ldots Q_n, G) := (\coprod_{i=1}^n Q_i,\rhd)$$ where $x \rhd y := x g_{i,j}(|y|)$ if $x \in Q_j, y \in Q_i$. We have the following theorem.

\begin{theorem} 
Let $(Q,\rhd)$ be a quandle. Then if it is not connected, $(Q,\rhd)$ may be expressed uniquely as $\#(Q_1,Q_2,\ldots, Q_n,
G)$ for some $G$, where $(Q_i, \rhd_i)$ are quandles, $Q_i$ are the orbits of the action of the inner automorphism group of $Q$ on $Q$, and $\rhd_i$ is the operation $\rhd$ restricted to $Q_i \times Q_i$.
\end{theorem}

\begin{proof}
Suppose $(Q, \rhd)$ is not connected. Since orbits of elements 
in $Q$ are the equivalence classes of an equivalence relation, $Q$ is uniquely expressed as a disjoint union of orbits in $Q$ under the group action of $Inn(Q)$, say $Q = \coprod_{i=1}^n Q_i$. Since each $Q_i$ is an orbit, we have
that $Q_i \rhd Q = Q_i$. Define $\rhd_i := \rhd|_{Q_i \times Q_i}$
($\rhd$ restricted to $Q_i \times Q_i$).  Then for each $i = 1
\ldots n$, $(Q_i,\rhd_i)$ is a quandle. ($\rhd_i$ inherits the
quandle structure of $\rhd$, so all that is still required to show that 
$(Q_i,\rhd_i)$ to be a quandle is closure.) 

Now we define the entries of matrix $G$, $g_{i,j}$.  For $Q_i$, we
have the augmentation map $|\cdot|_{Q_i}$ from $Q_i$ to
$\Gamma_{Q_i}$. For each $x \in Q_i$, let $g_{i,j}(|x|_{Q_i}) :=
\phi^j_x$, a right-action, where $\phi^j_x : Q_j \rightarrow Q_j$ and $y \phi^j_x =
y \rhd x$ for each $y \in Q_j$. Extend the map $g$ so that $g$ is a
homomorphism from $\Gamma_{Q_i}$ to $Aut(Q_j)$. (For $|x_1|_{Q_i},
|x_2|_{Q_i}, \ldots, |x_n|_{Q_i} \in \Gamma_{Q_1}$, let
$g_{i,j}(|x_1|^{\pm 1} |x_2|^{\pm 1} \ldots |x_n|^{\pm 1}) =
[g_{i,j}(|x_1|)]^{\pm 1}[g_{i,j}(|x_2|)]^{\pm 1} \ldots
[g_{i,j}(|x_n|)]^{\pm 1}$.) 

We must now check that $g_{i,j}$ is well-defined.  $g_{i,j}(|x \rhd
y|)$ is defined both as $\phi^j_{x \rhd y}$ and
$g_{i,j}(|y|^{-1}|x||y|) = [g_{i,j}(|y|)]^{-1} g_{i,j}(|x|)
g_{i,j}(|y|) = (\phi^j_y)^{-1} \phi^j_x \phi^j_y$.  Choose arbitrary
$z \in Q_j$.  Then $z \phi^j_y \phi^j_{x \rhd y} = (z \rhd y) \rhd
(x \rhd y) = (z \rhd x) \rhd y = z \phi^j_x \phi^j_y$. Since $z \in
Q_j$ was arbitrary, we see that $\phi^j_y \phi^j_{x \rhd y} 
\phi^j_x \phi^j_y$, i.e. $\phi^j_{x \rhd y} = (\phi^j_y)^{-1}
\phi^j_x \phi^j_y$. Hence $g_{i,j}$ is well-defined. \\ For all $y,
z \in Q_j$, $(y \rhd z)\phi^j_x = (y \rhd z) \rhd x = (y \rhd x)
\rhd (z \rhd x) = (y \phi^j_x) \rhd (z \phi^j_x)$. Also since $Q$ is
a quandle, we see that for each $q \in Q_j$ and $y \in Q_i$, there
exists a unique element, $q \mbox{ }\rhd^{-1} \mbox{ } y \in
Q_j$, such that $(q \mbox{ } \rhd^{-1} \mbox{ } y)
g_{i,j}(|y|) = q$.  Thus the image of $g_{i,j}$ is in $Aut(Q_j)$. 

So $Q = \#(Q_1,Q_2,\ldots, Q_n, G)$ by construction. 

Now we show that this decomposition is unique.  Earlier, we showed
that each $Q_i$ is uniquely determined.  Since $Q = \#(Q_1,Q_2,\ldots, Q_n,
G)$, we see that for each $x,y \in Q_i$, $x \rhd_i y$ must equal $x
\rhd y$.  Thus $\rhd_i$ is uniquely determined as well. Hence each subquandle 
$(Q_i, \rhd_i)$ is uniquely determined. For each $y \in Q_i$ and $x \in Q_j$,
we have that $x \rhd y = x g_{i,j}(|y|_{Q_i})$. Thus
$g_{i,j}(|y|_{Q_i})$ is uniquely determined for each $y \in Q_i$.
(An automorphism is determined by where it takes each element of its
domain to.) Since $g_{i,j}$ is a homomorphism on $\Gamma_{Q_i}$, and
$\Gamma_{Q_i}$ is generated by the elements $|q|_{Q_i}$ where $q \in
Q_i$, we have that $g_{i,j}$ is uniquely determined on
$\Gamma_{Q_i}$.  Hence this decomposition is unique up to
re-ordering.
\end{proof}

Now, unlike more familiar decomposition or factorization theorems, while
this decomposition is unique, it does not decompose the quandle into
indecomposable pieces, since, while they are single orbits under
$Inn(Q)$, the $Q_i$ may not be single orbits under their own groups
of inner automorphisms.  Nonetheless, iterating the construction of
the previous theorem, every quandle can be iteratively decomposed 
into connected quandles.  The uniqueness result of the previous theorem
then gives the uniqueness of the iterative decomposition.  

Of course, an arbitrary matrix $G$ of group homomorphisms $g_{i,j}:\Gamma_{Q_i}
\rightarrow Aut(Q_j)$ need not give rise to a
quandle.  We now give necessary and sufficient for $\#(Q_1,Q_2,\ldots Q_n,
G)$ to be a quandle:

\begin{theorem}\footnote{Shortly after the conclusion of the KSU REU, 
Nelson and Wong \cite{NW} announced
the independent discovery of a decomposition theorem
equivalent to this theorem and the preceding, when viewed as a
decomposition theorem, rather than a construction.   In their
result, the extra structure is expressed in terms of
compatible rack actions, rather than 
compatible group homomorphisms
from universal augmentation groups to automorphism groups.}
\label{meshthm}
Let $(Q_i,\rhd_i)$ be quandles for $i = 1 \ldots n$, and let
$g_{i,j}$ be homomorphisms from $\Gamma_{Q_i}$ to $Aut(Q_j)$ for 
$1 \leq i, j \leq n$, with $g_{i,i}$ the canonical homomorphism
from $\Gamma_{Q_i}$ to $Aut(Q_i)$ . Then \\ $\#(Q_1,Q_2,\ldots Q_n, G):=
(Q, \rhd)$ is a quandle if and only if for all $i, j, k$
distinct, $1 \leq i, j, k \leq n$, the following conditions hold:
\begin{enumerate}
\begin{item}
$(x g_{j,i}(|y|_{Q_j})) \rhd_i z = (x \rhd_i z) g_{j,i}(|y
g_{i,j}(|z|_{Q_i})|_{Q_j})$
\end{item}
\begin{item}
$(x g_{j,i}(|y|_{Q_j})) g_{k,i}(|z|_{Q_k}) = (x g_{k,i}(|z|_{Q_k}))
() g_{j,i}(|y g_{k,j}(|z|_{Q_k})|_{Q_j})$
\end{item}
\end{enumerate}
\end{theorem}

\begin{proof}
Clearly $Q$ is closed under the operation $\rhd$.  Since $Q_1, Q_2,
\ldots Q_n$ are quandles, we see that $x \rhd x = x$ for all $x \in
Q$.  Let $y, z \in Q$ be arbitrary.  If $y, z \in Q_i$ for some $i$,
then we see that $(z \mbox{ } \rhd^{-1} \mbox{ } y) \rhd y =
z$.  Alternatively, suppose $y \in Q_i$, $z \in Q_j$ for some $i
\neq j$. Then since $g_{i,j}(|y|_{Q_i})$ is an automorphism from
$Q_j$ to $Q_j$, there exists a unique $x \in Q_j$ such that $x
g_{i,j}(|y|_{Q_i}) = z$.  Thus $Q$ is a quandle if and only if the
third property holds: $(x \rhd y) \rhd z = (x \rhd z) \rhd (y \rhd
z)$ for all $x, y, z \in Q$.  If $x, y, z \in Q_i$ for some $Q_i$,
this is given.  Thus we have 4 cases:
\begin{enumerate}
\begin{item}
$x, y \in Q_i$, $z \in Q_j$
\end{item}
\begin{item}
$x, z \in Q_i$, $y \in Q_j$
\end{item}
\begin{item}
$y, z \in Q_i$, $x \in Q_j$
\end{item}
\begin{item}
$x\in Q_i$, $y \in Q_j$, $z \in Q_k$.
\end{item}
\end{enumerate}
\underline{\textbf{Case 1:}}  $(x \rhd y) \rhd z = (x \rhd_i y)
g_{j,i}(|z|_{Q_j})$. \\ $(x \rhd z) \rhd (y \rhd z) = (x
g_{j,i}(|z|_{Q_j})) \rhd_i (y g_{j,i}(|z|_{Q_j}))$. But these
coincide since $g_{j,i}(|z|)$ is an automorphism of $Q_i$. \\
\\
\underline{\textbf{Case 2:}}  $(x \rhd y) \rhd z = (x
g_{j,i}(|y|_{Q_j})) \rhd_i z$. \\ $(x \rhd z) \rhd (y \rhd z) = (x
\rhd_i z) g_{j,i}(|y g_{i,j}(|z|_{Q_i})|_{Q_j})$.  Thus for all $i
\neq j$, $1 \leq i, j \leq n$, $$(x g_{j,i}(|y|_{Q_j})) \rhd_i z =
(x \rhd_i z) g_{j,i}(|y g_{i,j}(|z|_{Q_i})|_{Q_j}),$$ which is condition (1). \\ \\
\underline{\textbf{Case 3:}}  $(x \rhd
y) \rhd z = (x g_{i,j}(|y|_{Q_i})) g_{i,j}(|z|_{Q_i})$.
\begin{align*}
(x \rhd z)\rhd (y \rhd z) & = (x g_{i,j}(|z|_{Q_i})) g_{i,j}(|y
\rhd_i z|_{Q_i})
\\ & = x g_{i,j}(|z|_{Q_i}|y \rhd_i z|_{Q_i}) 
\\ & \hspace*{1cm}\mbox{ since } g_{i,j}
\mbox{ is a group homomorphism.}
\\ & = x g_{i,j}(|z|_{Q_i}|z|_{Q_i}^{-1} |y|_{Q_i}
|z|_{Q_i}) 
\\ & \hspace*{1cm}\mbox{ since } |\cdot|_{Q_i} \mbox{ is an augmentation.}
\\ & = x g_{i,j}(|y|_{Q_i} |z|_{Q_i})
\\ & = (x g_{i,j}(|y|_{Q_i})) g_{i,j}(|z|_{Q_i}).
\end{align*}
Thus in this case, $(x \rhd y) \rhd z = (x \rhd z) \rhd (y \rhd z)$.
\\ \\
\underline{\textbf{Case 4:}}  $(x \rhd y) \rhd z = (x
g_{j,i}(|y|_{Q_j})) g_{k,i}(|z|_{Q_k})$. \\ $(x \rhd z) \rhd (y \rhd
z) = (x g_{k,i}(|z|_{Q_k})) g_{j,i}(|y g_{k,j}(|z|_{Q_k})|_{Q_j})$.
Thus for all $i, j, k$ distinct, $1 \leq i, j, k \leq n$, $$x
g_{j,i}(|y|_{Q_j})) g_{k,i}(|z|_{Q_k}) = (x g_{k,i}(|z|_{Q_k}))
g_{j,i}(|y g_{k,j}(|z|_{Q_k})|_{Q_j}),$$ which is condition (2). \\ \\

Therefore the given conditions are necessary and sufficient for $(Q, \rhd)$
to be a quandle.
\end{proof}

In principle, at least, this result reduces the classification problem
to the classification problem for connected quandles.  
To generate every quandle, we
iteratively specify which quandles to compose, and use the
conditions above to determine all matrices $G$ for which the
composition is a quandle.

The previous theorem also provides a means of constructing
new quandles from old:

\begin{definition}
A {\em mesh} for a sequence of quandles $Q_1,\ldots , Q_n$
is a matrix $G$ of group homomorphism $g_{ij}:\Gamma_{Q_i}\rightarrow Aut(Q_j)$ satisfying
the conditions of Theorem \ref{meshthm}.
\end{definition}

Notice that the matrix whose off-diagonal entries are
each the appropriate trivial group homomorphism, with
diagonal entries given by the canonical homomorphisms
is always a mesh.

Finally we name the construction described above:

\begin{definition}
Given a finite sequence of quandles $Q_1, \ldots , Q_n$
and a mesh $G$, the quandle $\#(Q_1,\ldots, Q_n, G)$
is the {\em semidisjoint union of $Q_1, \ldots , Q_n$
with respect to $G$}.

For $G$ the mesh with trivial off-diagonal entries,
$\#(Q_1,\ldots, Q_n, G)$ is the {\em disjoint union of
$Q_1, \ldots , Q_n$}.
\end{definition}

A final note before turning to connected
quandles:  the disjoint union of quandles is not the
coproduct in the category of quandles--the coproduct is the quotient of the quandle freely generated by the disjoint
union by the congruence which enforces all equations
holding in the individual quandles.  

\section{On the Classification of Connected Quandles}

Here we investigate conditions on a group for it to arise as the
group of inner automorphisms of a connected quandle, and derive a structure
theorem relating connected quandles and their groups of inner automorphisms.
In this section, we denote $Inn(Q)$ by $G$, where $Q$ is an connected 
quandle. Note that this $G$ is not related to the matrix of group homomorphisms in the previous section.

The key to the structure theorem is the fact that connected quandles are
single orbits of their inner automorphism groups, and thus by standard
results can be identified
as $G$-sets with a homogeneous space of cosets.  In particular

\begin{proposition}
Let $Q$ be an connected quandle on $n$ elements.  Then $n$ divides the order of $G$, and, moreover, any choice of  $q \in Q$ induces a $G$-equivariant bijection between $Q$ and $H \backslash G$, where $H$ is the stablizer of $q$.
\end{proposition}

In the case of $n$ prime, the converse of the first conclusion also holds:

\begin{theorem}
Let $Q$ be a quandle with $p$ elements, where $p$ is prime.  Then
$Q$ connected $\Leftrightarrow$ $p$ divides the order
of $G$.
\end{theorem}
\begin{proof}
$(\Rightarrow)$ This follows immediately from the above lemma. \\
$(\Leftarrow)$ Suppose $p$ divides $|G|$.  Then $|G| = p^a b$ for some positive
integers $a$ and $b$.  By the 1st Sylow theorem, $G$ has a subgroup
of order $p^a$. Choose one such subgroup $H$.  Hence every element in $H$ has order which
divides $p^a$.  Since $G$ and thus $H$ is a subgroup of ${\mathfrak S}_p$, every
element of $H$ has order which divides $p!$. Thus every element of
$H$ has order $1$ or $p$.  Choose an element of order $p$.
An element of $G \subset {\mathfrak S}_p$ which has order $p$ must be of the
form $(\begin{matrix} a_1 & a_2 & \ldots & a_p
\end{matrix})$ where the $a_i$s are the $p$ distinct elements of quandle $Q$.
In particular, $a_{i+1}$ is in the orbit of $a_i$ under the group
action of $G$. Hence $a_1, a_2, \ldots, a_n$ are in the same orbit
under the group action of $G$, i.e. $G$ is connected.
\end{proof}

By the Proposition, $Q$ is in bijection with the right cosets $H \backslash G$.
Hence we may represent $Q$ by $(\{Hg_1, Hg_2, Hg_3, \ldots, Hg_n\},
\rhd)$, where $Hg_i \rhd Hg_j = Hg_k$ if $q_i \rhd q_j = q_k$.

Assume $Q$ has the representation $(\{Hg_1, Hg_2, Hg_3, \ldots,
Hg_n\}, \rhd)$.  We define an augmentation map, $|\cdot|: Q = H
\backslash G \rightarrow G$, such that $Hg_i$ is mapped to $g \in G$
which takes $x \in Q$ to $x \rhd Hg_i$ for $i = 1, \ldots, n$.  To distinguish between this augmentation and existing notation for the order of a group, consider $|H|$ to be the order of the subgroup H and $|Hh|$ to be the augmentation of H as a right-coset in $H \backslash G$. Also denote the center of $H$ as $Z(H)$.

\begin{theorem} \label{gen}
Let $Q$ be an connected quandle on $n$ elements.  Let
$G = Inn(Q)$, $H$, $g_i$, and $|\cdot|$ be defined as above.  Then
$H \subset G \subset {\mathfrak S}_n$, $\frac{|G|}{|H|}= n$, $|Hh| \in Z(H)$, and $G$ is
generated by $|Hh|, |Hg_2|, \ldots,$ $|Hg_n|$, where $|H g_i| =
g_i^{-1} |Hh| g_i$.
\end{theorem}

\begin{proof}
By construction, $H \subset G$. $G$ is contained in the group
of bijective maps from the elements of $Q$ to the elements of $Q$, so $G \subset {\mathfrak S}_n$.  Since there are exactly $n$ cosets of $H$ in $G$,
we see that $\frac{|G|}{|H|}= n$.  Also by definition, $G = Inn(Q)$ is
generated by $|Hh|, |H g_2|, \ldots, |H g_n| \in G$.  It remains to
prove that $|Hh| \in Z(H)$ and $|H g_i| = g_i^{-1} |Hh| g_i$.

\begin{claim}
For all $g \in G$, $|H g_i g| = g^{-1} |H g_i| g$.
\end{claim}

\begin{proof}
The RHS maps $Hg_j g$ to $(Hg_j \rhd Hg_i) g$.  The LHS maps $Hg_j
g$ to $Hg_j g \rhd Hg_i g$.  But since $g \in Inn(Q)$, 

\small{
\begin{eqnarray*} (Hg_j \rhd
Hg_i) g & = & (\ldots (Hg_j \rhd Hg_i) \rhd^{\pm 1} Hg_{k_1} )\rhd^{\pm
1} Hg_{k_2}) \ldots \rhd^{\pm 1} Hg_{k_n}) \ldots ) \\
&  = & (\ldots(Hg_j
\rhd^{\pm 1} Hg_{k_1}) \rhd^{\pm 1} Hg_{k_2}) \ldots \rhd^{\pm 1}
Hg_{k_n})\ldots) \\
& & \hspace*{1cm} \rhd (\ldots(Hg_i \rhd^{\pm 1} Hg_{k_1}) \rhd^{\pm
1} Hg_{k_2}) \ldots \rhd^{\pm 1} Hg_{k_n})\ldots) 
\\ &= & Hg_j g \rhd Hg_i g.
\end{eqnarray*}}

\noindent  Hence the LHS and the RHS coincide.
\end{proof}

Choosing $g$ to be the identity element in $G$, we see that $|H g_i|
= g_i^{-1} |Hh| g_i$.  Taking $i = 1$ and $g \in H$ arbitrary, we see
that $|Hh| = |H g| = |H g_1 g| = g^{-1} |H g_1| g = g^{-1} |Hh| g$.
Hence $|Hh|$ commutes with any element in $H$. Since $H = H g_1 = H
g_1 |H g_1| = H |Hh|$, we see that $|Hh| \in H$. Hence $|Hh| \in Z(H)$.
\end{proof}

\begin{theorem} \label{structure}
Suppose that for groups $G$ and $H$, we have that $H \subset G
\subset {\mathfrak S}_n$, $\frac{|G|}{|H|}= n$.  Let $g_1, g_2, \ldots, g_n$ be coset
representatives of $H$ in $G$.  Suppose also that $G$ is generated
by $g_1^{-1}|Hh|g_1$, $g_2^{-1}|Hh|g_2$, $g_3^{-1}|Hh|g_3$, \ldots,
$g_n^{-1}|Hh|g_n$ for some $|Hh| \in Z(H)$. Then $Hg_i \rhd Hg_j =
Hg_i g_j^{-1} |Hh| g_j$ defines an connected quandle
with $n$ elements.
\end{theorem}

\begin{proof}
Note that if $H g_j = H g'_j$, then $g_j^{-1} |Hh| g_j = (h
g'_j)^{-1} |Hh| (h g'_j) = (g'_j)^{-1} h^{-1} |Hh| h g'_j$ for some $h
\in H$.  But since $|Hh| \in Z(H)$, this is equal to $(g'_j)^{-1} |Hh|
g'_j$.  Hence the staitment that $G$ is generated by $g_1^{-1} |Hh|$,
$g_2^{-1}|Hh|g_2$, $g_3^{-1}|Hh|g_3$, \ldots, $g_n^{-1}|Hh|g_n$ for
some $|Hh| \in Z(H)$ makes sense, and also $Hg_i \rhd Hg_j$ is
well-defined.  We now check to see if this defines a quandle.  $H
g_i \rhd H g_i = H g_i g_i^{-1} |Hh| g_i = H |Hh| g_i = H g_i$ since
$|Hh| \in Z(H)$.  Now for $j, k$ arbitrary, $H g_i \rhd H g_j = H
g_k$ implies that $H g_i = H g_k g_j^{-1} |Hh|^{-1} g_j$.  Note that
such an $i$ exists and is unique.  Finally, we have that $(H g_i
\rhd H g_k) \rhd (H g_j \rhd H g_k) = (H g_i g_k^{-1} |Hh| g_k) \rhd
(H g_j g_k^{-1} |Hh| g_k)$.  \\
Let $H g_j g_k^{-1} |Hh| g_k = H g_m$
for some $m$.  Then $g_m = h g_j g_k^{-1} |Hh| g_k$ for some $h \in
H$. Hence

\begin{align*}
(H g_i \rhd H g_k) \rhd (H g_j \rhd H g_k) & = (H g_i g_k^{-1} |Hh|
g_k) \rhd H g_m 
\\ & = H g_i g_k^{-1} |Hh| g_k g_m^{-1} |Hh| g_m 
\\ &
= H g_i g_k^{-1} |Hh| g_k g_k^{-1} |Hh|^{-1} g_k g_j^{-1}
\\ & \hspace*{1cm} h^{-1} |Hh| h
g_j g_k^{-1} |Hh| g_k 
\\ & = H g_i g_j^{-1} h^{-1} |Hh| h g_j g_k^{-1}
|Hh| g_k 
\\ & = H g_i g_j^{-1} |Hh| g_j g_k^{-1} |Hh| g_k 
\\ & = (H g_i \rhd
H g_j) \rhd H g_k.
\end{align*}

Since $g_1^{-1}|Hh|g_1$, $g_2^{-1}|Hh|g_2$, $g_3^{-1}|Hh|g_3$, \ldots,
$g_n^{-1}|Hh|g_n$ generate $G$, we see that for all $i, j$, $g_i^{-1}
g_j \in G$ is generated by $g_1^{-1}|Hh|g_1$, $g_2^{-1}|Hh|g_2$,
$g_3^{-1}|Hh|g_3$, \ldots, $g_n^{-1}|Hh|g_n$.  Hence for each $i, j$,
there exist $k_1, k_2, \ldots, k_n$ such that $$H g_j = H g_i
g_i^{-1} g_j = (\ldots (H g_i \rhd H g_{k_1}) \rhd H g_{k_2}) \rhd
\ldots \rhd H g_{k_n},$$ i.e. the quandle defined by $Hg_i \rhd Hg_j
= Hg_i g_j^{-1} |Hh| g_j$ is connected.
\end{proof}

The above theorems provide a program for constructing
all finite connected quandles more satisfactory than the brute
force approach of \cite{HN}: for each $n$, test the generation condition
of
Theorem \ref{structure} for all triples of a subgroup $G$ of
${\mathfrak S}_n$, a 
subgroup $H \subset G$, and a central element of $H$.

We conclude with a number of restrictions on groups that 
arise as inner automorphism groups of finite quandles, which follow
easily from Theorem \ref{structure}:

\begin{corollary}
For $(Q,\rhd)$ connected and $n > 1$, $G$ is not abelian.
\end{corollary}

\begin{proof}  We proceed by contradiction:
suppose $G$ were abelian.  Then for all $g \in G$, 
$|Hg| = |H g_1 g| =
g^{-1} |Hg_1| g = |Hg_1| = |Hh|$.  Since $|H g_i|$ fixes $Hg_i$ for
$i = 1, \ldots, n$, $|Hh|$ must also fix $Hg_i$ for each $i = 1,
\ldots, n$.  Hence $|Hh|$ is trivial, and thus $|Hg_i|$ is trivial
for all $i$. But, since the $G$ is generated by the $|H g_i|$,$G$ is trivial.  
Therefore $(Q,\rhd)$ is not connected unless $n =1$, which contradicts 
the hypothesis that $n > 1$. Thus, the corollary holds.
\end{proof}

\begin{corollary}
For $Q$ connected on $n$ elements, $G = {\mathfrak S}_n$ implies
that $n = 1, 3$.
\end{corollary}

\begin{proof}
Since $\frac{|G|}{|H|}= n$, then $|H| = (n-1)!$.  But $H \subset {\mathfrak S}_{n-1}$,
so $H = {\mathfrak S}_{n-1}$.  For $n \geq 4, Z({\mathfrak S}_{n-1}) = {\bf 1}$, the trivial one-element group.  Since $|Hh| \in Z(H)$, $|Hh|$ must equal ${\bf 1}$.  But then $|H g_i| = {\bf 1}$ for all $i$, so
$G = {\bf 1}$, a contradiction.  Hence $n \leq 3$.  For $n = 2$, $H = {\bf 1}$
and $G = {\mathfrak S}_2 \neq {\bf 1}$, a contradiction as above.  This gives the
desired result.
\end{proof}

\begin{corollary}
For $Q$ connected on $n$ elements and $n > 1$, $G =
A_n$ implies that $n = 4$.
\end{corollary}

\begin{proof}
Since $\frac{|G|}{|H|} = n$, then $|H| = \frac{1}{2}(n-1)!$.  But $H \subset
{\mathfrak S}_{n-1}$ and $H \subset A_n$, so $H = A_{n-1}$.  For $n \geq 5$,
$Z(A_{n-1}) = {\bf 1}$. Since $|Hh| \in Z(H)$, $|Hh|$ must equal ${\bf 1}$.  But
then $|H g_i| = 1$ for all $i$, so $G = 1$, a contradiction.  Hence
$n \leq 4$.  But, for $n = 2, 3$, $A_n$ of order $1$ or is abelian. 
\end{proof}

The results of this section suggest an approach to classifying
finite connected quandles which should be more computationally effective
than the brute force approach of Ho and Nelson \cite{HN}: for
a given order, determine (up to conjugacy) all towers of groups $H \subset G
\subset {\mathfrak S}_n$ for which the index of $H$ in $G$
is $n$, and in  which $H$ has a non-trivial center.  Central
elements of $H$ can then be tested against the generating
condition of Theorems \ref{gen} and \ref{structure}.

\begin{flushleft}
\bibliography{Book}
\end{flushleft}

\end{document}